\documentclass{l4dc2025}

\makeatletter
\renewcommand{\jmlrvolume}[1]{}
\renewcommand{\jmlryear}[1]{}
\renewcommand{\jmlrproceedings}[2]{}
\def\ps@jmlrtps{%
	\def\@oddhead{}
	\def\@evenhead{}
	\def\@oddfoot{}
	\def\@evenfoot{}
}
\makeatother

\usepackage{wrapfig}
\usepackage{hyperref}
\usepackage{cleveref}
\usepackage{graphicx}
\usepackage{enumitem}
\usepackage{ifthen}
\usepackage{graphicx}
\usepackage{float}
\usepackage{caption}
\usepackage{wrapfig}
\usepackage[compact]{titlesec} 
\usepackage{bm}
\usepackage{tikz}
\usepackage{algorithm,algorithmic}
\usepackage{comment}
\usepackage{tikz-cd}
\usepackage{comment}
\newtheorem{Theorem}{Theorem}
\newtheorem{Lemma}{Lemma}
\newtheorem{Definition}{Definition}
\newtheorem{Corollary}{Corollary}
\newtheorem{Assumption}{Assumption}
\theorembodyfont{\upshape}
\newtheorem{Remark}{Remark}
\DeclareMathOperator*{\argmin}{arg\,min}

\title[Adaptive Control of Positive Systems with Application to Learning SSP]{\hspace{-1mm}Adaptive Control of Positive Systems with Application to Learning SSP}
\usepackage{times}
\author{%
 \Name{Fethi Bencherki} \Email{fethi.bencherki@control.lth.se} \\
 \Name{Anders Rantzer} \Email{anders.rantzer@control.lth.se} \\
 \addr  Department of Automatic Control,
Lund University, Sweden
}

\begin{document}

\maketitle

\begin{abstract}
An adaptive controller is proposed and analyzed for the class of infinite-horizon optimal control problems in positive linear systems presented in \citep{ohlin2024optimal}. This controller is derived from the solution of a ``data-driven algebraic equation" constructed using the model-free Bellman equation from Q-learning. The equation is driven by data correlation matrices that do not scale with the number of data points, enabling efficient online implementation. Consequently, a sufficient condition guaranteeing stability and robustness to unmodeled dynamics is established. The derived results also provide a quantitative characterization of the interplay between excitation levels and robustness to unmodeled dynamics. The class of optimal control problems considered here is equivalent to Stochastic Shortest Path (SSP) problems, allowing for a performance comparison between the proposed adaptive policy and model-free algorithms for learning the stochastic shortest path, as demonstrated in the numerical experiment.
\end{abstract}
\section{Introduction}
Positive systems represent a class of dynamical systems in which the state remains nonnegative for all time, provided the initial condition is nonnegative. Many physical variables—such as concentrations, buffer levels, queue lengths, charge levels, and prices—are inherently nonnegative, motivating mathematical models that incorporate this structural constraint. Such models have been successfully employed in diverse application domains, including traffic and congestion modeling \citep{shorten2006positive}, thermodynamics \citep{haddad2010nonnegative}, biology and pharmacology \citep{carson2013modelling,blanchini2014piecewise}, and epidemiology \citep{hernandez2013modeling}. These wide-ranging applications have spurred significant research interest in the analysis and control of positive systems. Foundational and recent contributions in this area include \citep{de2001stabilization,rami2007controller,rantzer2015kalman, ebihara2016analysis,gurpegui2023minimax,ohlin2024optimal}, as well as the comprehensive tutorial \citep{rantzer2018tutorial}.

Despite the extensive literature on the optimal control of positive systems, most existing approaches remain offline and model-based. This work is motivated by recent advancements and growing interest in statistical machine learning \citep{tsiamis2023statistical} and finite-time analysis, aiming to develop an online, data-driven approach to the problem. This direction has gained increasing attention in recent years, as seen in studies addressing the linear quadratic regulator (LQR) problem \citep{de2019formulas,markovsky2021behavioral, zhao2024data} and more recent works focused on positive systems \citep{shafai2022data,miller2023data,padoan2023data,iwata2024data,makdah2023model,wang2024data2}.

In the current manuscript, we adopt a worst-case approach regarding the types of disturbances and uncertain plant parameters, consistent with the works presented in \citep{rantzer2021minimax,kjellqvist2022learning,kjellqvist2022minimax,bencherki2023robust,renganathan2023online}. However, the disturbances in this paper are subject to a bounded constraint based on past states and inputs, similar to the approach in \citep{rantzer2024data}, where the online linear quadratic optimal control problem in the presence of non-stochastic process disturbances is addressed. This constraint introduces a different perspective compared to the aforementioned works.

\subsection{Contributions and outline of the paper}
\paragraph{Contributions}
The paper proposes and analyzes an adaptive control scheme for the class of optimal control of positive systems presented in \citep{ohlin2024optimal}. A data-driven algebraic equation is constructed from the model-free counterpart of the Bellman optimality equation. This equation allows for the direct extraction of the adaptive policy, bypassing an explicit system identification step. The data-driven equation is constructed and updated under the assumption of available full-state measurements contaminated with additive process noise. To account for uncertainties related to unmodeled dynamics and time-variations in the plant, we avoid making stochastic assumptions about the noise. This ensures that the proposed approach is equipped with robustness guarantees.

\paragraph{Applications}
The considered problem class can capture various network routing problem settings. A specific instance of this appeared in \citep{bencherki2024data}, where the authors addressed the problem of learning the optimal processing rate over processing networks. Another interesting instance is the Stochastic Shortest Path (SSP) problem class, due to the work in \citep{ohlin2024heuristic}, where the authors demonstrate the equivalence between the two problem classes. This equivalence motivates a numerical study comparing the performance of the model-free policy presented here with existing parameter-free algorithms for SSP.

\paragraph{Outline}
The paper is structured as follows: Section 2 outlines the problem setup, starting with the model-based optimal control problem from \citep{ohlin2024optimal}, progressing to its model-free counterpart, and deriving the algebraic equation for the adaptive policy. Section 3 presents supporting lemmas, leading to the formal online performance analysis of the proposed adaptive policy. Section 4 compares the numerical performance of the adaptive policy with the $Q$-learning algorithm from \citep{yu2013boundedness} in learning the stochastic shortest path. Finally, Section 5 concludes the paper and discusses future directions. 

\subsection{Notation}
Inequalities are applied element-wise to matrices and vectors throughout. Furthermore, the notation $\mathbb{R}^n_+$ denotes the closed positive orthant of dimension $n$. The symbol $\lvert X \rvert$ represents the element-wise absolute value of the entries of the matrix (or vector) $X$. The operator $\min\{A,0\}$ extracts the minimum element of $A$, yielding zero if $A$ contains no negative elements. The expressions $\mathbf{1}_{p\times q}$ and $\mathbf{0}_{p\times q}$ represent matrices of ones or zeros, respectively, of the indicated dimension, with the subscript omitted when the size is clear from the context. The operation $\otimes$ denotes the Kronecker product. $U(a,b)$ denotes the uniform distribution between $a$ and $b$.

\section{Problem setup}
We consider the class of infinite-horizon optimal control problems presented in \citep{ohlin2024optimal}:
{\small
\begin{equation}
\label{pb1}
\begin{aligned}
    \textnormal{Minimize} \quad & \sum_{t=0}^{\infty} \left[ s^\top x(t) + r^\top u(t) \right] \quad \textnormal{over } \{u(t)\}_{t=0}^{\infty} \\
    \textnormal{subject to} \quad & x(t+1) = Ax(t) + Bu(t) \\
    & u(t) \ge 0, \quad x(0) = x_0 \\
    & 
    \begin{matrix}
        \mathbf{1}^\top u_1(t) & \le & E_1^\top x(t) \\
        \vdots & & \vdots \\
        \mathbf{1}^\top u_n(t) & \le & E_n^\top x(t)
    \end{matrix}
\end{aligned}
\end{equation}}
where  \
$A \in \mathbb{R}^{n \times n}$ and $B = \begin{bmatrix} B_1 & \cdots & B_n \end{bmatrix} \in \mathbb{R}^{n \times m}$, where each $B_i \in \mathbb{R}^{n \times m_i}$, define the linear dynamics. The input signal $u \in \mathbb{R}^m$ is partitioned into $n$ subvectors $u_i$, each containing $m_i$ elements, such that $m = \sum_{i = 1}^{n} m_i$. The cost vectors associated with the states and control inputs are $s \in \mathbb{R}^{n}_{+}$ and $r \in \mathbb{R}^{m}_{+}$, where each $r_i \in \mathbb{R}^{m_i}_{+}$ follows the partitioning of $u$. The constraints on the input signal $u$ are given by $E = \begin{bmatrix} E_1 & \cdots & E_n \end{bmatrix}^\top \in \mathbb{R}^{n \times n}_{+}$. Furthermore, we define the extended constraint matrix 
$
\bar{E} = \begin{bmatrix}
    {\bf 1}_{m_1}^\top \otimes E_1 \quad \cdots \quad {\bf 1}_{m_n}^\top \otimes E_n
\end{bmatrix}^\top \in \mathbb{R}_+^{m \times n},
$
and the set of indices 
$
\mathcal{V} = \{1, \ldots, n\}.
$
 Let  
$
K = \begin{bmatrix} K_1^\top & \cdots & K_n^\top \end{bmatrix}^\top
$
be a feedback matrix with $K_i \in \mathbb{R}^{m_i \times n}_{+}$, and define the set of feasible gains as
\begin{equation}\label{gains}
    \mathcal{K}(E) \triangleq  \left\{ K : (\forall i \in \mathcal{V})\; \mathbf{1}_{m_i}^\top K_i = E_i^\top \;\textnormal{or}\; \mathbf{1}_{m_i}^\top K_i = \mathbf{0}_{1\times n} \right\}.
\end{equation}
Correspondingly, the state feedback law for the \( i \)-th control subvector is given by \( u_i = K_i x \). The set \( \mathcal{K}(E) \) characterizes all feedback gain matrices that result in either full or zero actuation of the control inputs \( u_i \), for \( i \in \mathcal{V} \).
We impose the following two assumptions on the sextuple \( (A, B, E, \bar{E}, s, r) \).
\begin{Assumption}[\citep{ohlin2024heuristic}]
\label{as:ABE}
    The matrices $A,B$ and the set $\mathcal{K}\left(E\right)$ satisfy 
    $
        (A+BK)x \ge 0
    $
     for all $K\in\mathcal{K}\left(E\right)$ and all reachable states $x\in\mathbb{R}^n_+$.
\end{Assumption}
\begin{Assumption}\label{as:sr}
The triplet \( (s, \bar{E}, r) \) satisfies \( s > \bar{E}^\top r \).
\end{Assumption}
\begin{Remark}
Assumption~\ref{as:ABE} ensures the positivity of the closed-loop dynamics. Assumption~\ref{as:sr}, on the other hand, requires that \( s > 0 \), which guarantees detectability of the system and ensures that the optimal feedback law is also stabilizing \citep{ohlin2024heuristic}.
\end{Remark}
\subsection{Solution to problem~\ref{pb1}}
Under Assumption~\ref{as:ABE} and via dynamic programming \citep{bellman1966dynamic}, it was shown in~\citep[Theorem~1]{ohlin2024optimal} that if problem (\ref{pb1}) has a finite value for every \( x(0) \in \mathbb{R}^n_+ \), then the optimal cost would be \( p^\top x(0) \), where \( p \in \mathbb{R}_+^n \) is the solution to the following model-based algebraic equation
\begin{align}
        \label{ARE}
            p = s + A^\top p + \sum_{i = 1}^n \min \{ r_i +B_i^\top p, 0 \} E_i.
\end{align}
Furthermore, the optimal policy is a linear state feedback law \( u(t) = Kx(t) \), where 
\( K = \begin{bmatrix} K_{1}^\top ~\ldots~ K_{n}^\top \end{bmatrix}^\top \) and  
\begin{align} \label{opt-gain}
    K_i \triangleq \begin{bmatrix} \mathbf{0}_{j-1\times n} \\ E_i^\top \\ \mathbf{0}_{m_i-j\times n} \end{bmatrix} \in  \mathbb{R}^{m_i \times n},
\end{align}
where the vector \( E_i^\top \) enters at the \( j \)-th row, with \( j \) being the index of the minimal element of \( r_i + B_i^\top p \), provided it is negative. If all elements are nonnegative, then \( K_i = \mathbf{0}_{m_i\times n} \). The solution to (\ref{ARE}) can, for instance, be obtained via value iteration (VI), by performing the following fixed-point iteration on the parameter \( p \) until convergence
\begin{align}
    \label{value_itr_p}
    p^{k+1} = s + A^\top p^{k} + \sum_{i = 1}^M \min \{ r_i + B_i^\top p^k, 0 \} E_i, \quad p^0 = 0.
\end{align}
\subsection{Model-free optimal control of positive systems via $Q$-factor}
The cost-to-go function from time \( t \) for the optimization problem in (\ref{pb1}) under a control policy \( u \), starting at time \( t \) from state \( x(t) \), is given by  
\begin{align*}
   J( x(t))\triangleq \min_u \sum_{k=t}^{\infty}  {s}^\top x(k) + r^\top u(k).
\end{align*}
If finite, this optimization problem yields the objective value  
$
J(x(t)) = p^\top x(t),
$
where \( p \) is determined by solving (\ref{ARE}). The optimal \( Q \)-function, as defined in \citep{bradtke1994adaptive}, is  
\begin{align}\label{Q-def}
    Q(x(t),\, u(t)) \triangleq c(x(t),\, u(t)) + J(x(t+1)),
\end{align}
which represents the cost of taking action \( u(t) \) starting at state \( x(t) \) and subsequently following the optimal policy \( u^* \). The optimal $Q$-function is then given by
 \begin{align}\label{Q-fct2}
    Q^*( x(t),\, u(t))=\begin{bmatrix}
s^\top +   p^\top A & r^\top+  p^\top B
\end{bmatrix}\begin{bmatrix}
 x(t) \\
  u(t)
\end{bmatrix} = \begin{bmatrix}
q^x\\
 q^u
\end{bmatrix}^\top \begin{bmatrix}
 x(t) \\
  u(t)
\end{bmatrix}
=q^\top\begin{bmatrix}
 x(t) \\
  u(t)
\end{bmatrix},
\end{align}
 where
 $
    q\triangleq\begin{bmatrix}
q^x\\
 q^u
\end{bmatrix}\triangleq\begin{bmatrix}
s  +   A^\top  p \\ r+ B^\top  p
\end{bmatrix}
$. It also holds that  
$
    J(x(t)) = \min_{ u(t) \in \mathcal{U}(x(t))} Q( x(t),\,u(t)),
$
where \( \mathcal{U}(x(t)) \) denotes the set of inputs satisfying the constraints in (\ref{pb1}) at time \( t \).
 Then, from (\ref{Q-def}), it follows  
\begin{align}\label{bel-Q}
    Q( x(t),\, u(t)) = c(x(t),\,u(t)) + \min_{ u(t+1) \in \mathcal{U}(x(t+1))} Q( x(t+1),\,u(t+1)),  
\end{align}  
and the optimal policy is given by  
$
    u^*(t) = \argmin_{ u(t) \in \mathcal{U}(x(t))} Q( x(t),\,u(t))
$. One can interpret (\ref{bel-Q}) as the Bellman equation in the \( Q \)-factor formulation \citep{sutton2018reinforcement}.  
By virtue of the definition in (\ref{gains}), replacing (\ref{Q-fct2}) in (\ref{bel-Q}) yields
\begin{align}\label{Q-learn1}
    \begin{bmatrix}
q^x \\ q^u
\end{bmatrix}^\top\begin{bmatrix}
 x(t) \\
 u(t)
\end{bmatrix}=s^\top x(t)+r^\top u(t)+ \min_{ K\in \mathcal{K}\left(E\right)}  \begin{bmatrix}
q^x \\ q^u
\end{bmatrix}^\top\begin{bmatrix}
I \\
 K
\end{bmatrix}x(t+1).
\end{align}
In the absence of knowledge of the dynamics \( (A, B) \), equation (\ref{Q-learn1}) enables us to obtain information about the \( q \)-parameter by collecting triplets \( (x(t),\, u(t),\, x(t+1)) \). In fact, collecting \( t+1 \) consecutive data points, for any \( t \geq 1 \), leads to
\begin{align}\label{Ric-x3}
    \left(  {q}-\begin{bmatrix}
        s\\ r
    \end{bmatrix}\right)^\top
    \begin{bmatrix}
        x(0) & \cdots &  x(t-1)\\
        u(0) & \cdots &  u(t-1)
    \end{bmatrix} = \min_{ K\in \mathcal{K}(E)} q^\top  
    \begin{bmatrix}
        I\\ K
    \end{bmatrix}     
    \begin{bmatrix}
        x(1) & \cdots &  x(t)
    \end{bmatrix}.
\end{align}
Multiplying equation~\eqref{Ric-x3} from the right by  
$
\begin{bmatrix}
    \lambda^{t-1} x(0) & \lambda^{t-2} x(1) & \cdots & x(t-1) \\
    \lambda^{t-1} u(0) & \lambda^{t-2} u(1) & \cdots & u(t-1)
\end{bmatrix} \in \mathbb{R}^{(n+m) \times t}
$
for a forgetting factor \( \lambda \in (0,\,1] \), and defining the data correlation matrices
\begin{align}
    \Sigma(t) \triangleq \sum_{k=0}^{t-1} \lambda^{t-1-k} 
    \begin{bmatrix}
        x(k) \\
        u(k)
    \end{bmatrix}
    \begin{bmatrix}
        x(k) \\
        u(k)
    \end{bmatrix}^\top + \lambda^t \Sigma(0) \;\; \text{and} \;\; 
    \bar{\Sigma}(t) \triangleq \sum_{k=0}^{t-1} \lambda^{t-1-k} x(k+1)
    \begin{bmatrix}
        x(k) \\
        u(k)
    \end{bmatrix}^\top, \label{Sigma}
\end{align}
yields the \emph{data-driven algebraic equation} in the \( q(t) \)-parameter
\begin{align}\label{DDARE}
    \left( q(t) - 
    \begin{bmatrix}
        s \\ r
    \end{bmatrix} \right)^\top \Sigma(t) =
    \min_{K(t) \in \mathcal{K}(E)} q(t)^\top 
    \begin{bmatrix}
        I \\ K(t)
    \end{bmatrix} \bar{\Sigma}(t).
\end{align}
Here, we denote the data-based solution by \( (q(t), K(t)) \), in contrast to the model-based (or ground-truth) solution \( (q, K) \). Equation~\eqref{DDARE} forms the foundation for the construction of the proposed controller, as will be shown in the sequel.
\subsection{Problem Formulation}
Inspired by (\ref{DDARE}), for the linear system 
\begin{align}\label{noisy-state}
    x(t+1) = A x(t) + B u(t) + w(t),
\end{align}
we propose and analyze the performance of policies of the form 
\begin{align}\label{cert-equiv}
    \begin{cases}
        \Sigma(t) = \lambda \Sigma(t-1) + 
        \begin{bmatrix}
            x^\top(t-1) & u^\top(t-1)
        \end{bmatrix}^\top
        \begin{bmatrix}
            x^\top(t-1) & u^\top(t-1)
        \end{bmatrix}, \quad \Sigma(0) \succ 0,\\
        \bar{\Sigma}(t) = \lambda \bar{\Sigma}(t-1) + x(t) 
        \begin{bmatrix}
            x^\top(t-1) & u^\top(t-1)
        \end{bmatrix}, \quad \bar{\Sigma}(0) = 0,\\
        u(t) = K(t) x(t) + \epsilon(t).
    \end{cases}
\end{align}
The controller states \( \Sigma(t) \) and \( \bar{\Sigma}(t) \) accumulate correlation data using a forgetting factor \( \lambda \). Based on these statistics, the control gain \( K(t) \) is computed as the minimizing argument of the data-driven algebraic equation in~\eqref{DDARE}. To ensure sufficient excitation for learning the true system dynamics \( (A, B) \), additive exploration noise \( \epsilon(t) \) is introduced.
\begin{Definition}\label{Def1}
  Let the model parameter set \( \mathcal{M}_\beta \) be the set of plants \( (A, B, E) \) satisfying Assumptions~\ref{as:ABE} and~\ref{as:sr}, such that the algebraic equation in (\ref{ARE}) admits a solution \( p \) satisfying  
\begin{align*}
    s \leq p \leq (\beta \min_{i} s_i  ) {\mathbf 1} \leq \beta s.
\end{align*}
\end{Definition}
\begin{Remark}
    The parameter \( \beta \) reflects the degree of stabilizability of the system. Larger values of \( \beta \) correspond to systems that are harder to stabilize.
\end{Remark}
Constrained to the model set defined in Definition~\ref{Def1}, we aim to establish guarantees for the stability and robustness of the closed-loop system under the policy given in~\eqref{cert-equiv}, assuming certain properties of the disturbance \( w \). Before proceeding with the analysis, we first present methods for solving the algebraic data-driven equation in~\eqref{DDARE}, as its solution is central to both the construction and implementation of the proposed controller.
\subsection{Solution to data-driven algebraic equation in the $q(t)$-parameter}
Similarly to the \( p \)-parameter-based algebraic equation in~\eqref{ARE}, the equation in~\eqref{DDARE} can be solved using value iteration, policy iteration, or linear programming. For brevity, policy iteration is omitted.
\paragraph{Solution via value iteration} 
The solution is obtained by performing the following fixed-point iteration on the \( q(t) \)-parameter, starting from \( q^0(t) = 0 \), and continuing until convergence:
{\small \begin{align}
    q^{k+1}(t) &= \Sigma^{-1}(t) \bar\Sigma^\top(t) 
    \begin{bmatrix}
        I & \left( K^k(t) \right)^\top 
    \end{bmatrix} q^k(t) 
    + 
    \begin{bmatrix}
        s \\ r
    \end{bmatrix}, 
    \quad q^0(t) = 0, \label{value-itr1} \\[-1ex]
    K^k_i(t) &= 
    \begin{cases}
        \begin{bmatrix}
            \mathbf{0}_{j-1 \times n} \\ 
            E_i^\top \\ 
            \mathbf{0}_{m_k - j \times n}
        \end{bmatrix}, & \text{if } \min\left\{ \left(q_i^u\right)^k(t),\, 0 \right\} < 0 \\
        \mathbf{0}_{m_i \times n}, & \text{otherwise}
    \end{cases}, \quad i = 1, \ldots, n, \nonumber
\end{align}}
where \( j \) denotes the index of the most negative entry in \( \left(q_i^u\right)^k(t) \) and the full controller matrix is
$
K^k(t) = 
\begin{bmatrix} 
\left(K_{1}^k(t)\right)^\top ~\ldots~ \left(K_{n}^k(t)\right)^\top 
\end{bmatrix}^\top.
$
\paragraph{Solution via linear programming}
Instead of using a fixed-point iteration, the solution can also be obtained via a linear programming (LP) formulation. To this end, define  
$
z \triangleq \begin{bmatrix} (z^x)^\top & (z^u)^\top \end{bmatrix}^\top,
$
where \( z^x \in \mathbb{R}_+^{n} \) and \( z^u \in \mathbb{R}_+^{m} \). We propose the following optimization problem to solve for \( q(t) \):
{\small \begin{align*}
    \text{maximize} \quad & \mathbf{1}^\top 
    \begin{bmatrix}
        I & -E^\top
    \end{bmatrix} z \quad \text{over } z \in \mathbb{R}_+^{n+m},\; q(t) \in \mathbb{R}^{n+m} \\
    \text{subject to} \quad & \Sigma^{-1}(t)\, \bar\Sigma^\top(t) 
    \begin{bmatrix}
        I & -E^\top
    \end{bmatrix} z = q(t) - 
    \begin{bmatrix}
        s \\ r
    \end{bmatrix}, \\
    & z^x = q^x(t), \quad z^u \geq q^u(t), \quad z^u \geq 0.
\end{align*}}
\begin{Remark}
    When the plant \( (A, B) \) is known, the term \( \Sigma^{-1}(t)\, \bar\Sigma^\top(t) \) is replaced by the model \( \begin{bmatrix} A & B \end{bmatrix}^\top \). A corresponding model-based version of the fixed-point iteration in~\eqref{value-itr1} becomes:
   {\small \begin{align}
        q^{k+1} &=  
        \begin{bmatrix}
            A & B
        \end{bmatrix}^\top 
        \begin{bmatrix}
            I & \left( K^k \right)^\top 
        \end{bmatrix} q^k 
        + 
        \begin{bmatrix}
            s \\ r
        \end{bmatrix}, 
        \quad q^0 = 0, \label{value-itr-q2} \\[-1ex]
        K^k_i &= 
        \begin{cases}
            \begin{bmatrix}
                \mathbf{0}_{j-1 \times n} \\ 
                E_i^\top \\ 
                \mathbf{0}_{m_k - j \times n}
            \end{bmatrix}, & \text{if } \min\left\{ \left(q_i^{u}\right)^k,\, 0 \right\} < 0 \\
            \mathbf{0}_{m_i \times n}, & \text{otherwise}
        \end{cases}, \quad i = 1, \ldots, n, \nonumber
    \end{align}}
    where \( K^k = \begin{bmatrix} \left(K_1^k\right)^\top ~ \ldots ~ \left(K_n^k\right)^\top \end{bmatrix}^\top \), and \( j \) is the index of the most negative element in \( \left(q_i^u\right)^k \). Here, the subscript \( t \) is omitted from \( q \) because the time-dependent and optimal \( q \)-parameters coincide.
\end{Remark}

\section{Main results}
We begin by establishing an important supporting result.
\begin{Lemma}\label{lemma1}
    Iterating on \( q \) in~\eqref{value-itr-q2} is algebraically equivalent to iterating on \( p \) in~\eqref{value_itr_p}.
\end{Lemma}
\textbf{Proof}. See Appendix-\ref{lem1}. \jmlrQED
\begin{Remark} \label{rem4}
We extract two key observations from Lemma~\ref{lemma1}:
\begin{enumerate}[label=(\roman*)]
    \item Lemma~\ref{lemma1} asserts that value iteration in the \(q\)-parameter is algebraically equivalent to value iteration in the \(p\)-parameter. \label{Rem1a}
    \item According to Lemma~\ref{lemma1}, given \( (q(t), K(t)) \), the solution to the data-driven algebraic equation in~\eqref{DDARE}, if we define
$
p(t) \triangleq \begin{bmatrix} I & \left(K(t)\right)^\top \end{bmatrix} q(t),
$
then \( p(t) \) satisfies the following model-based algebraic equation under the estimated dynamics \( (\hat{A}(t), \hat{B}(t)) \), provided that persistently exciting data is collected. Specifically, defining
$
[\hat{A}(t) \;\; \hat{B}(t)] \triangleq \bar{\Sigma}(t) \Sigma^{-1}(t),
$
we obtain
\begin{align}\label{data-based}
    p(t) - s = \hat{A}^\top(t) p(t) + \sum_{i = 1}^n \min \{ r_i + \hat{B}_i^\top (t) p(t), 0 \} E_i.
\end{align}
\end{enumerate}
This is key in establishing proofs for the main results of the paper, as shall be seen in Section~\ref{PA}.
\end{Remark}
\begin{Lemma}\label{lemma3}
 Let \( p \in \mathbb{R}_+^n \) and \( \hat{p} \in \mathbb{R}_+^n \) be such that they satisfy the following algebraic equation and inequality, respectively, i.e., they satisfy
    \begin{align*}
                        p = s + A^\top p + \sum_{i = 1}^n \min \{ r_i +B_i^\top p, 0 \} E_i \quad \text{and} \quad
                         \hat p \ge s + A^\top \hat p + \sum_{i = 1}^n \min \{ r_i +B_i^\top \hat p, 0 \} E_i.
    \end{align*}
 Then, it holds that $\hat p \ge p$. 
\end{Lemma}
\textbf{Proof}. See Appendix-\ref{lem2}. \jmlrQED
\subsection{Performance analysis}\label{PA}
In anticipation of our first main result, we define
\begin{align*}
    \tilde{\Sigma}(t) \triangleq 
    \begin{bmatrix} 
        \Sigma^{wx}(t) & \Sigma^{wu}(t) 
    \end{bmatrix} 
    \triangleq \sum_{k=0}^{t-1} \lambda^{t-1-k} w(k) 
    \begin{bmatrix} 
        x^\top(k) & u^\top(k) 
    \end{bmatrix},
\end{align*}
which implies that 
$
\bar{\Sigma}(t) = 
\begin{bmatrix} 
    A & B 
\end{bmatrix} \Sigma(t) + \tilde{\Sigma}(t)
$. Therefore,
\begin{align}\label{Def-Dyn}
    \begin{bmatrix} 
        \hat{A}(t) & \hat{B}(t) 
    \end{bmatrix} = \bar{\Sigma}(t) \Sigma^{-1}(t)
    = \begin{bmatrix} 
        A & B 
    \end{bmatrix} + 
    \begin{bmatrix} 
        \tilde{A}(t) & \tilde{B}(t) 
    \end{bmatrix},
\end{align}
where
$
\begin{bmatrix} 
    \tilde{A}(t) & \tilde{B}(t) 
\end{bmatrix} \triangleq \tilde{\Sigma}(t) \Sigma^{-1}(t) 
= \begin{bmatrix} 
    \Sigma^{wx}(t) \Sigma^{-1}(t) & \Sigma^{wu}(t) \Sigma^{-1}(t) 
\end{bmatrix}
$
represents the model misspecification. The first main result of the paper is stated next.
\begin{Theorem}\label{thm1}
    Consider \( \beta,\,\rho \in \mathbb{R}_+ \) satisfying \( \rho \beta < 1 \), and let \( \Sigma(t) \) and \( \bar{\Sigma}(t) \) be as in (\ref{Sigma}). Let \( \left(A,\,B\right) \in \mathcal{M}_\beta \). Additionally, suppose that 
     \begin{align}\label{pert}
         \left\lVert E^\top \right\rVert_\infty \left\lVert A^\top - \begin{bmatrix} I & 0 \end{bmatrix} \Sigma^{-1}(t) \bar{\Sigma}^\top(t) \right\rVert_\infty + \left\lVert B^\top - \begin{bmatrix} 0 & I \end{bmatrix} \Sigma^{-1}(t) \bar{\Sigma}^\top(t) \right\rVert_\infty \leq \rho, \;\; \forall t \ge t_0.
     \end{align}
     Let \( p \in \mathbb{R}_+^n \) be the solution to the model-based algebraic equation in (\ref{ARE}), and let \( q(t) \) be the solution of the data-based algebraic equation in (\ref{DDARE}) with \( K(t) \) being the minimizing argument. Define \( p(t) \triangleq \begin{bmatrix} I & K^\top(t) \end{bmatrix} q(t) \). Then, it holds that 
     \begin{align}\label{two-sides}
         \hat{\alpha} p \leq p(t) \leq \check{\alpha}^{-1} p
     \end{align}
     for positive constants satisfying \( \check{\alpha} = 1 - \rho \beta \) and \( \hat{\alpha} = 1 - \check{\alpha}^{-1} \rho \beta \).
\end{Theorem}
\textbf{Proof}. See Appendix-\ref{thm1-pf}. \jmlrQED
\begin{Remark}
	The result in theorem~\ref{thm1} is obtained from a perturbation analysis to the solution of the algebraic equations in (\ref{ARE}) and (\ref{data-based}). Perturbation analysis of algebraic Riccati equations is a well studied problem in the literature, see the works \citep{konstantinov1993perturbation,sun1998perturbation,konstantinov2003perturbation}.
\end{Remark}
	Assumption~(\ref{pert}) can be expressed as  
{\small $
	\left\lVert E^\top \right\rVert_\infty  \left\lVert \Sigma^{-1}(t) \left({\Sigma}^{wx}(t)\right)^\top\right\rVert_\infty + \left\lVert \Sigma^{-1}(t) \left({\Sigma}^{wu}(t)\right)^\top\right\rVert_\infty \leq \rho
$}
and holds as long as the following two conditions are satisfied:  
\begin{enumerate}[label=(\roman*)]
	\item The condition number of \( \Sigma(t) \) must not be too large, which depends on the choice of a suitable exploration signal \( \epsilon(t) \) and a sufficiently large \( t_0 \) to ensure sufficient data is collected.
	\item The disturbance sequence \( w(t) \) remains sufficiently small relative to \( x(t) \). If \( w(t) \) represents unmodeled dynamics, this implies that such dynamics are subject to a gain bound from the state to the disturbance.
\end{enumerate}
In cases where the disturbance sequence \( w \) is modeled as i.i.d.\ Gaussian noise, increasing \( t_0 \) typically enables Assumption~(\ref{pert}) to hold for smaller values of \( \rho \). However, if the system is exposed to adversarial disturbances, the true dynamics may never be accurately learned, as the adversary \( w \) could adopt a policy that persistently misleads the controller.

\begin{Theorem}\label{thm2}
Consider \( \beta,\,\rho \ge 0 \) satisfying \( \rho \beta < 1 \), and let \( \Sigma(t) \) and \( \bar{\Sigma}(t) \) be as defined in~\eqref{Sigma}. Let \( (A, B) \in \mathcal{M}_\beta \), and let \( p \in \mathbb{R}_+^n \) denote the solution to the model-based algebraic equation in~\eqref{ARE}. Suppose that~\eqref{pert} holds for all \( t \ge t_0 \), and let \( K(t) \) denote the minimizer in~\eqref{DDARE}. Then, it holds that
\begin{align}\label{theta-ineq}
    \check{\alpha}^{-1} \left(1 + \rho \beta \left( 1 + \beta \left\lVert A + \left\lvert B \right\rvert \bar{E} \right\rVert_1 \right) \right) p - s 
    \ge A^\top p + K^\top(t) \left( r + B^\top p \right),
\end{align}
where the constant \( \check{\alpha} \) is as defined in Theorem~\ref{thm1}.
\end{Theorem}
\textbf{Proof}. See Appendix-\ref{Thm2}. \jmlrQED
\begin{Remark}
Theorem~\ref{thm2} aims to characterize how the storage function \( p^\top x \) is affected when the data-driven gain \( K(t) \) is used in place of the optimal gain \( K \). As \( \rho \to 0 \), we have \( \check{\alpha} \to 1 \), and the right-hand side of (\ref{theta-ineq}) approaches 1, thereby closing the suboptimality gap. Larger values of \( \beta \) correspond to systems that are more difficult to stabilize, necessitating smaller values of \( \rho \) to satisfy the condition \( \rho \beta < 1 \), which in turn reflects the need for higher-quality data.

\end{Remark}
\begin{Corollary}\label{coro1}
Consider \( \left(A,\,B\right) \in \mathcal{M}_\beta \), with \( p \) being the positive solution to the algebraic equation~(\ref{ARE}), and let \( \beta,\,\rho > 0 \) satisfy \( \rho \beta < 1 \). Let the system evolve as
\begin{align*}
    x(t+1) = A x(t) + B u(t) + w(t), \quad\quad u(t) = K(t) x(t) + \epsilon(t),
\end{align*}
where \( \Sigma(t) \) and \( \bar{\Sigma}(t) \) are defined as in~(\ref{Sigma}), and assume that~(\ref{pert}) holds for all \( t \ge t_0 \). Let \( q(t) \) be the solution to the data-based algebraic equation in~(\ref{DDARE}), with \( K(t) \) being its minimizing argument. Then, the following inequality holds:
\begin{align*}
    \sum_{t=t_0}^{T-1} s^\top x(t) + r^\top K(t) x(t) 
    \leq \gamma^{-1} \left( p^\top x_{t_0} + \sum_{t=t_0}^{T-1} \beta\, s^\top \left\lvert B \epsilon(t) + w(t) \right\rvert \right),
\end{align*}
where \( 0 < \gamma \leq 1 \) satisfies
$
    \gamma \left( s - \bar{E}^\top \lvert r \rvert \right) 
    \leq s - \bar{E}^\top \lvert r \rvert 
    - \left( \check{\alpha}^{-1} - 1 
    + \rho \beta \left( 1 + \beta \left\lVert A + \left\lvert B \right\rvert \bar{E} \right\rVert_1 \right) \right) \beta s,
$
and the constant \( \check{\alpha} \) is as defined in Theorem~\ref{thm1}.
\end{Corollary}
\textbf{Proof}. See Appendix-\ref{Coro1}. \jmlrQED
\begin{Remark}
  Corollary~\ref{coro1} quantifies the suboptimality in the incurred cost when applying the adaptive controller in place of the optimal one. The constant \( \gamma \) reflects the impact of model misspecification and satisfies \( \gamma \to 1 \) as \( \rho \to 0 \), thereby tightening the suboptimality bound. The final term on the right-hand side accounts for the accumulated cost due to external disturbances \( w(t) \) and the exploration noise \( \epsilon(t) \).
\end{Remark}
\section{Numerical experiment}
We consider the problem of learning the following stochastic shortest path
\begin{center}
{\small 
\textbf{Stochastic Shortest Path (SSP) problem} 
\begin{align*}
                \mathcal{T}\left(1\right)=\begin{bmatrix}
            0.4 & 0 \\
            0 & 0.4 \\
            0.4 & 0.4  \\
            0.2 & 0.2 
        \end{bmatrix},\,
        \mathcal{T}\left(2\right)=\begin{bmatrix}
            0 & 0.3 & 0\\
            0.6 & 0 & 0.1\\
            0.4 & 0.7 & 0.4 \\
            0 & 0 & 0.5
        \end{bmatrix},
        \mathcal{T}\left(3\right)=\begin{bmatrix}
            0 & 0.2 \\
            0 & 0.2 \\
            0.4 & 0 \\
            0.6 & 0.6
        \end{bmatrix},\,        \mathcal{T}\left(4\right)=\begin{bmatrix}
            0  \\
            0  \\
            0 \\
            1 
        \end{bmatrix}\\
        c\left(1\right)=\begin{bmatrix}
            1.5 & 2
        \end{bmatrix},\,c(2)=\begin{bmatrix}
            1.5 & 2 & 2
        \end{bmatrix},\,c(3)=\begin{bmatrix}
            1.5 & 2
        \end{bmatrix} c(4)=0,\, i_{\text{init}}=1
\end{align*}}
\end{center}
\begin{center}
\hspace{25mm}
  \begin{tikzpicture}
\draw[color=blue,line width=1mm,-{Triangle[length=3mm, width=3mm]}](-0.3,0) -- ++ (0,-1.5);
 \node[text width=4cm] at (1.9,-0.7) {\textcolor{blue}{Conversion according to {\citep{ohlin2024heuristic}} }};
     \end{tikzpicture}  
\end{center}
\begin{center}
\vspace{-2mm}
{\textbf {Reformulation as problem~\ref{pb1}}}\\ 
\vspace{0.2cm}
 $n = 3$, $m = 4$, $x_0 = \left[1 \; 0 \; 0\right]^\top$ and dynamics
    \begin{align*}
        A = \begin{bmatrix}
            0.4 & 0 & 0 \\ 0 & 0.6 & 0 \\ 0.4 & 0.4 & 0.4
        \end{bmatrix} , \,
        B = \begin{bmatrix}
            -0.4 & 0.3 & 0 & 0.2 \\ 0.4 & -0.6 & -0.5 & 0.2 \\ 0 & 0.3 & 0 & -0.4
        \end{bmatrix}
    \end{align*}
    \textnormal{with associated costs}
    \begin{align*}
        s &= \begin{bmatrix}
            1.5 & 1.5 & 1.5
        \end{bmatrix}^\top, \;\;\;
        r = \begin{bmatrix}
            0.5 & 0.5 & 0.5 & 0.5
        \end{bmatrix}^\top\nonumber\\
        E&=I,\,m_1=1,\,m_2=2,\,m_3=1
    \end{align*}
\end{center}
where \( \mathcal{T}(i) \), for \( i \in \mathcal{V} \cup \{v_g\} \) with \( v_g = 4 \) denoting the fictitious goal state, represent the transition probabilities from state \( i \) to other states, where different columns correspond to different actions. The vector \( c(i) \) denotes the expected stage cost associated with transitioning from state \( i \) to other states under different actions. 
Note that both Assumptions~\ref{as:ABE} and~\ref{as:sr} hold true for our system. We compare the performance of the adaptive policy to that of the \( Q \)-learning algorithm presented in~\citep{yu2013boundedness}, employing an \(\epsilon\)-decreasing exploration strategy with \( \epsilon = 0.05 \alpha^h \), where \( \alpha = 0.99 \) and \( h \) denotes the episode number.
In the context of Problem~\ref{pb1}, \(\epsilon\)-greedy exploration means selecting a random gain \( K \in \mathcal{K}(E) \) with probability \( \epsilon \), and choosing the estimated optimal gain \( K(t) \), computed via the data-driven algebraic equation~(\ref{DDARE}), with probability \( 1 - \epsilon \).
For the policy in~(\ref{cert-equiv}), we set \( \lambda = 1 \) and initialize \( \Sigma(0) = 10^{-6}I \). To make the learning task more challenging, we corrupt the state measurements with additive positive disturbances \( w(t) \sim U(0,\, 0.01) \) for all \( t \).
The comparison between the two algorithms is performed by evaluating the regret, defined as
\begin{align*}
    R(H) = \left(\sum_{h=0}^{H-1}\sum_{t=0}^{T_h-1}s^\top x(t) + r^\top u(t)\right) - HJ(x(0)),
\end{align*}
where \( H \) denotes the number of episodes, and \( T_h \) is the duration of episode \( h \). The term \( J(x(0)) \) 

\begin{wrapfigure}[25]{r}{0.5\textwidth}
    \centering
    \subfigure[{\small Accumulated regret of each algorithm with $\epsilon$-decreasing explorations where $\epsilon=0.05 \alpha^h$ for $\alpha=0.99$ and $h$ the episode number.} ]
        {\label{fig1}  \includegraphics[width=0.9\linewidth]{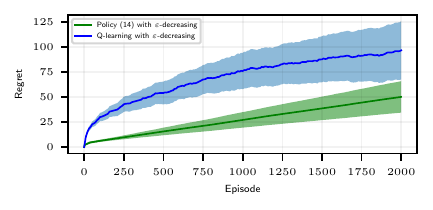}}
    \subfigure[{\small The condition in~(\ref{pert}) evaluated for \( \rho = 0.3 \) in the presence of disturbances. For this example, \( \lVert E^\top \rVert_\infty = 1 \).
} ]{\label{fig2}
          \includegraphics[width=0.9\linewidth]{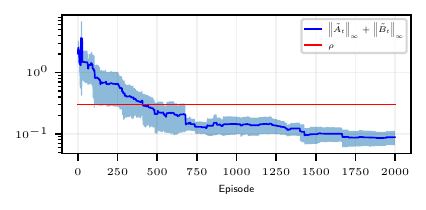}
    }
    \caption{{\small Each plot represents the average over 100 repeated runs, with the shaded area indicating the 95\% confidence interval.
}}
\end{wrapfigure}
 \noindent represents the optimal cumulative cost obtained by applying the optimal policy to the system subject to disturbances \( w \).
In the SSP domain, an episode terminates when the goal state (\( i = 4 \)) is reached, whereas in Problem~\ref{pb1}, an episode ends when the state vector becomes sufficiently small, indicating that the measurements are dominated by noise. After each episode, the system states are re-initialized in both problem domains. The algorithms are then restarted, retaining the latest estimates: the \( Q \)-factor for the \( Q \)-learning algorithm, and the \( q(t) \)-parameter and \( \Sigma(t) \) for the policy in~(\ref{cert-equiv}). For more details on the implementation, refer to the code\footnote{Code available at: \href{https://github.com/Fethi-Bencherki/adaptive-control-positive-systems-l4dc2025}{https://github.com/Fethi-Bencherki/adaptive-control-positive-systems-l4dc2025}}.
The resulting regrets are shown in~\Cref{fig1}. Both algorithms exhibit sublinear regret, with our proposed algorithm consistently outperforming the \( Q \)-learning approach. \Cref{fig2} further illustrates a test of the condition in~(\ref{pert}), using \( \rho = 0.3 \). Notably, the condition is satisfied after a sufficient number of episodes, suggesting that more accurate estimates of the system are obtained as additional data is gathered.


\section{Conclusion}
The paper presented a robust adaptive data-driven control framework tailored towards the class of positive systems presented in \citep{ohlin2024optimal}. This was achieved via the construction of a data-driven algebraic equation in the $Q$-factor, based on which the controller policy is updated in an online fashion. The designed policy proved to robustly stabilize the set $\mathcal{M}_\beta$ with robustness meant to be tolerance to a certain degree of unmodeled dynamics. The considered class witnesses applications in network routing problems among which are Stochastic Shortest Path problems, allowing for performance comparison with the existing model-free methods of finding the stochastic shortest path. Future work concerns
exploring better exploration strategies than $\epsilon$-greedy, and the possibility of adapting efficient exploration methods from the SSP literature into our control setup.


\acks{
The authors thank David Ohlin at Lund University for his insightful feedback, which helped improve the paper. They are affiliated with the ELLIIT Strategic Research Area at Lund University. This project received funding from the European Research Council (ERC) under grant agreement No.~834142 (ScalableControl), and was partially supported by the Wallenberg AI, Autonomous Systems and Software Program (WASP), funded by the Knut and Alice Wallenberg Foundation.
}
\bibliography{refs}
\clearpage

\appendix
\section{Proofs}
\subsection{Proof of Lemma~\ref{lemma1}}\label{lem1}
Consider the iteration in~(\ref{value-itr-q2}) and define
\[
p^k \triangleq \begin{bmatrix}
I & \left(K^k(t)\right)^\top
\end{bmatrix} q^k.
\]
This then gives
\begin{align}\label{eq1}
q^{k+1} = \begin{bmatrix}
A^\top \\ B^\top
\end{bmatrix} p^k + \begin{bmatrix}
s \\ r
\end{bmatrix},
\end{align}
and also
\begin{align}\label{eq3}
p^{k+1} \triangleq \begin{bmatrix}
I & \left(K^{k+1}(t)\right)^\top
\end{bmatrix} q^{k+1}.
\end{align}
where the feedback matrix \( K^{k+1}_i \) is given by
\begin{align}\label{eq2}
K^{k+1}_i = \begin{cases}
\begin{bmatrix}
\mathbf{0}_{j-1 \times n} \\
E_i^\top \\
\mathbf{0}_{m_k - j \times n}
\end{bmatrix}, & \text{if } \min\left\{ \left(q_i^u\right)^{k+1},\,0 \right\} \stackrel{(\ref{eq1})}{=} \min\left\{ r_i + B_i^\top p^{k+1},\,0 \right\} < 0, \\
0_{m_i \times n}, & \text{otherwise}
\end{cases}, \quad i = 1, \ldots, n.
\end{align}
Plugging (\ref{eq1}) and (\ref{eq2}) into (\ref{eq3}) yields the desired result:
\begin{align*}
    p^{k+1} = s + A^\top p^k + \sum_{i = 1}^n \min \{ r_i + B_i^\top p^k,\, 0 \} E_i.
\end{align*}
\jmlrQED
\subsection{Proof of lemma~\ref{lemma3}}\label{lem2}
Start by defining the value functions \( J(x) \triangleq p^\top x \) and \( \hat{J}(x) \triangleq  \hat{p}^\top x \), which solve the Bellman equation and inequality, respectively. That is, they satisfy
\begin{align*}
    J(x) &= \min_{u \in \mathcal{U}(x)} \left[ s^\top x + r^\top u + J(Ax + Bu) \right], \\[-1ex]
    \hat{J}(x) &\ge \min_{\hat{u} \in \mathcal{U}(x)} \left[ s^\top x + r^\top \hat{u} + \hat{J}(Ax + B\hat{u}) \right],
\end{align*}
for all \( x \in \mathbb{R}_+^n \).
 Taking the minimizers \( u = \begin{bmatrix}
        u_1^\top & \ldots & u_n^\top
    \end{bmatrix}^\top \) with \( u_i = -K_i x \), and \( \hat{u} = \begin{bmatrix}
        \hat{u}_1^\top & \ldots & \hat{u}_n^\top
    \end{bmatrix}^\top \) with \( \hat{u}_i = -\hat{K}_i x \) for \( i = 1, \ldots, n \), where \( K_i \) is defined as in~\eqref{opt-gain} and \( \hat{K}_i \) is defined analogously, leads to
\begin{align*}
    \left( p - s - A^\top p - \sum_{i = 1}^n \min \left\{ r_i + B_i^\top p, 0 \right\} E_i \right)^\top x &= 0, \\
    \left( \hat{p} - s - A^\top \hat{p} + \sum_{i = 1}^n \min \left\{ r_i + B_i^\top \hat{p}, 0 \right\} E_i \right)^\top x &\ge 0.
\end{align*}
  These last equation and inequality are guaranteed to hold if the equation and inequality in the statement of the lemma hold and \( x \in \mathbb{R}_+^n \). It then follows that \( \hat{J}(x) \ge J(x) \) for all \( x \in \mathbb{R}_+^n \); see Proposition~1.(a) in~\citep{li2024exact}. Consequently, the elementwise ordering \( \hat{p} \ge p \) holds.{\flushright \jmlrQED}
\subsection{Proof of Theorem~\ref{thm1}}\label{thm1-pf}
First, and according to Remark~\ref{rem4}, we associate the solution of the noisy (respectively, noiseless) data-driven algebraic equation in~(\ref{DDARE}) with the following model-based algebraic equations in the estimated (respectively, true) model:
\begin{align}
    p(t) - s &=  \hat{A}^\top(t)\, p(t) + \sum_{i = 1}^n \min \{ r_i + \hat{B}_i^\top(t)\, p(t), 0 \} E_i,\label{datab} \\
    p - s &= A^\top p + \sum_{i = 1}^n \min \{ r_i + B_i^\top p, 0 \} E_i.\label{modelb}
\end{align}
The two equations coincide when the uncertainty vanishes, i.e., when \( \rho = 0 \).
\begin{enumerate}[label=(\roman*)]
      \item Consider the algebraic equation in~(\ref{modelb}). Substituting
\begin{align*}
    \begin{bmatrix}
        A & B
    \end{bmatrix}
    = 
    \begin{bmatrix}
        \hat{A}(t) & \hat{B}(t)
    \end{bmatrix}
    - 
    \begin{bmatrix}
        \tilde{A}(t) & \tilde{B}(t)
    \end{bmatrix}
\end{align*}
into~(\ref{modelb}) yields
\begin{align*}
    p - s 
    &= \left( \hat{A}(t) - \tilde{A}(t) \right)^\top p 
    + \sum_{i = 1}^n \min \left\{ r_i + \left( \hat{B}_i(t) - \tilde{B}_i(t) \right)^\top p, 0 \right\} E_i \\
    &= \left( \hat{A}(t) - \tilde{A}(t) \right)^\top p 
    + \sum_{i = 1}^n \min \left\{ r_i + \hat{B}_i^\top(t) p - \tilde{B}_i^\top(t) p, 0 \right\} E_i.
\end{align*}
Since
\[
\min \left\{ r_i + \hat{B}_i^\top(t)p - \tilde{B}_i^\top(t)p,\, 0 \right\} 
\ge \min \left\{ r_i + \hat{B}_i^\top(t)p,\, 0 \right\} 
- \left\lVert \tilde{B}_i^\top(t)p \right\rVert_\infty,
\]
using the fact that the $\infty$-norm is sub-multiplicative and applying the triangle inequality yields
\begin{multline*}
    p - s 
    \ge \hat{A}^\top(t)p 
    + \sum_{i=1}^n \min \left\{ r_i + \hat{B}_i^\top(t)p,\, 0 \right\} E_i 
    - \sum_{i=1}^n \left\lVert \tilde{B}_i^\top(t)p \right\rVert_\infty E_i 
    - \left\lvert \tilde{A}^\top(t)p \right\rvert \\
    \ge \hat{A}^\top(t)p 
    + \sum_{i=1}^n \min \left\{ r_i + \hat{B}_i^\top(t)p,\, 0 \right\} E_i 
    - \lVert p \rVert_\infty \sum_{i=1}^n \left\lVert \tilde{B}_i^\top(t) \right\rVert_\infty E_i 
    - \left\lvert \tilde{A}^\top(t) \right\rvert \lvert p \rvert.
\end{multline*}
Next, we have
\[
\sum_{i = 1}^n \left\lVert \tilde{B}_i^\top(t) \right\rVert_\infty E_i 
\leq E^\top \mathbf{1} \max_i \left\lVert \tilde{B}_i^\top(t) \right\rVert_\infty 
\leq \left\lVert E^\top \right\rVert_\infty \left\lVert \tilde{B}^\top(t) \right\rVert_\infty \mathbf{1},
\]
and
\[
\left\lVert \left\lvert \tilde{A}^\top(t) \right\rvert \lvert p \rvert \right\rVert_\infty 
\leq \left\lVert \tilde{A}^\top(t) \right\rVert_\infty \lVert p \rVert_\infty,
\]
which implies
\begin{multline*}
    \lVert p \rVert_\infty \sum_{i=1}^n \left\lVert \tilde{B}_i^\top(t) \right\rVert_\infty E_i 
    + \left\lvert \tilde{A}^\top(t) \right\rvert \lvert p \rvert 
    \leq \left( \left\lVert E^\top \right\rVert_\infty \left\lVert \tilde{B}^\top(t) \right\rVert_\infty 
    + \left\lVert \tilde{A}^\top(t) \right\rVert_\infty \right) \lVert p \rVert_\infty \mathbf{1} \\
    \leq \left( \left\lVert E^\top \right\rVert_\infty \left\lVert \tilde{B}^\top(t) \right\rVert_\infty 
    + \left\lVert \tilde{A}^\top(t) \right\rVert_\infty \right) \beta \min_{i} s_i \mathbf{1} \\
    \leq \left( \left\lVert E^\top \right\rVert_\infty \left\lVert \tilde{B}^\top(t) \right\rVert_\infty 
    + \left\lVert \tilde{A}^\top(t) \right\rVert_\infty \right) \beta s,
\end{multline*}
In consequence, it follows from~(\ref{pert}) and by defining 
\[
\check{\alpha} \triangleq 1 - \rho \beta
\]
that
\begin{align*}
    p - \check{\alpha}s 
    \ge \hat{A}^\top(t)p 
    + \sum_{i=1}^n \min \left\{ r_i + \hat{B}_i^\top(t)p,\, 0 \right\} E_i.
\end{align*}
Dividing both sides by $\check{\alpha}$ and using the fact that $r_i \ge 0$ yields
\begin{multline*}
    \check{\alpha}^{-1} p - s 
    \ge \hat{A}^\top(t)\check{\alpha}^{-1}p 
    + \sum_{i=1}^n \min \left\{ \check{\alpha}^{-1}\left(r_i + \hat{B}_i^\top(t)p\right),\, 0 \right\} E_i \\
    \ge \hat{A}^\top(t)\check{\alpha}^{-1}p 
    + \sum_{i=1}^n \min \left\{ r_i + \hat{B}_i^\top(t)\check{\alpha}^{-1}p,\, 0 \right\} E_i.
\end{multline*}

An application of Lemma~\ref{lemma3} with $\hat{p} = \check{\alpha}^{-1}p$ then reveals that
\begin{align}\label{ineq1}
    p(t) \le \check{\alpha}^{-1} p.
\end{align}
\item  Now we consider the algebraic equation in~(\ref{datab}). Plugging in
\[
\begin{bmatrix}
\hat{A}(t) & \hat{B}(t)
\end{bmatrix} 
= \begin{bmatrix}
A & B
\end{bmatrix} 
+ \begin{bmatrix}
\tilde{A}(t) & \tilde{B}(t)
\end{bmatrix}
\]
yields
\begin{multline*}
    p(t) - s 
    = \left( A + \tilde{A}(t) \right)^\top p(t) 
    + \sum_{i = 1}^n \min \left\{ r_i + \left( B_i + \tilde{B}_i(t) \right)^\top p(t),\, 0 \right\} E_i \\
    = \left( A + \tilde{A}(t) \right)^\top p(t) 
    + \sum_{i = 1}^n \min \left\{ r_i + B_i^\top p(t) + \tilde{B}_i^\top(t) p(t),\, 0 \right\} E_i \\
    \ge A^\top p(t) 
    + \sum_{i = 1}^n \min \left\{ r_i + B_i^\top p(t),\, 0 \right\} E_i 
    - \lVert p(t) \rVert_\infty \sum_{i = 1}^n \left\lVert \tilde{B}_i^\top(t) \right\rVert_\infty E_i 
    - \left\lvert \tilde{A}^\top(t) \right\rvert \lvert p(t) \rvert \\
    \ge A^\top p(t) 
    + \sum_{i = 1}^n \min \left\{ r_i + B_i^\top p(t),\, 0 \right\} E_i 
    - \left( \left\lVert E^\top \right\rVert_\infty \left\lVert \tilde{B}^\top(t) \right\rVert_\infty 
    + \left\lVert \tilde{A}^\top(t) \right\rVert_\infty \right) \lVert p(t) \rVert_\infty \mathbf{1}.
\end{multline*}
Employing~(\ref{ineq1}) and~(\ref{pert}), we obtain
\begin{align*}
    p(t) - \left(1 - \rho \beta \check{\alpha}^{-1}\right)s 
    \ge A^\top p(t) + \sum_{i = 1}^n \min \left\{ r_i + B_i^\top p(t),\, 0 \right\} E_i.
\end{align*}
Define $\hat{\alpha} \triangleq 1 - \rho \beta \check{\alpha}^{-1}$, and dividing both sides by $\hat{\alpha}$ gives
\begin{align*}
    \hat{\alpha}^{-1} p(t) - s 
    \ge A^\top \hat{\alpha}^{-1} p(t) + \sum_{i = 1}^n \min \left\{ r_i + B_i^\top \hat{\alpha}^{-1} p(t),\, 0 \right\} E_i.
\end{align*}
Applying Lemma~\ref{lemma3} with $\hat{p} = \hat{\alpha}^{-1} p(t)$ yields
\begin{align*}
    p \le \hat{\alpha}^{-1} p(t).
\end{align*}
\end{enumerate}
Consequently, we obtain the bounds
\begin{align*}
    \hat{\alpha} p \le p(t) \le \check{\alpha}^{-1} p,
\end{align*}
for positive constants $\check{\alpha} = 1 - \rho \beta$ and $\hat{\alpha} = 1 - \rho \beta \check{\alpha}^{-1}$.
\jmlrQED
\subsection{Proof of Theorem~\ref{thm2}}\label{Thm2} 
 The following holds by letting $\hat{B}(t) = B + \tilde{B}(t)$ and $\hat{A}(t) = A + \tilde{A}(t)$ in~(\ref{data-based}), and using $K^\top(t)$ as defined in~(\ref{value-itr1}):
\begin{multline*}
    p(t) - s 
    = \left( A + \tilde{A}(t) \right)^\top p(t) 
    + K^\top(t) \left( r + \left( B + \tilde{B}(t) \right)^\top p(t) \right) \\
    = A^\top p(t) + K^\top(t) \left( r + B^\top p(t) \right) 
    + \tilde{A}^\top(t) p(t) + K^\top(t) \tilde{B}^\top(t) p(t) \\
    = A^\top p + K^\top(t) \left( r + B^\top p \right) 
    + \left( \tilde{A}^\top(t) + K^\top(t) \tilde{B}^\top(t) \right) p(t) 
    + \left( A^\top + K^\top(t) B^\top \right) \left( p(t) - p \right) \\
    \ge A^\top p + K^\top(t) \left( r + B^\top p \right) 
    + \left( \tilde{A}^\top(t) + K^\top(t) \tilde{B}^\top(t) \right) p(t) 
    - \left( \left\lvert A^\top \right\rvert + \bar{E}^\top \left\lvert B^\top \right\rvert \right) \left\lvert p(t) - p \right\rvert,
\end{multline*}
since $\left\lvert K(t) \right\rvert \le \bar{E}$ and 
\[
A^\top + K^\top(t) B^\top \le A^\top + \bar{E}^\top \left\lvert B^\top \right\rvert.
\]
From~(\ref{two-sides}), we have
\[
\left\lvert  p(t) - p \right\rvert \leq \left(1 - \hat{\alpha} \right) p,
\]
since $\check{\alpha}^{-1} - 1 = 1 - \hat{\alpha}$. Then, we obtain
\begin{multline*}
    \left(A^\top + \bar{E}^\top \left\lvert B^\top \right\rvert \right) \left\lvert p(t) - p \right\rvert 
    \leq \left(1 - \hat{\alpha} \right) \left\lVert A^\top + \bar{E}^\top \left\lvert B^\top \right\rvert \right\rVert_\infty \lVert p \rVert_\infty {\bf 1} \\
    \leq \left(1 - \hat{\alpha} \right) \left\lVert A^\top + \bar{E}^\top \left\lvert B^\top \right\rvert \right\rVert_\infty \beta s 
    \leq \left(1 - \hat{\alpha} \right) \left\lVert A + \left\lvert B \right\rvert \bar{E} \right\rVert_1 \beta p.
\end{multline*}
Moreover, using~(\ref{pert}) we get
\[
\left(\tilde{A}^\top(t) + K^\top(t)\tilde{B}^\top(t) \right)p(t) 
\leq \rho \check{\alpha}^{-1} \beta s 
\leq \rho \beta \check{\alpha}^{-1} p,
\]
which yields
\[
p(t) - s \ge A^\top p + K^\top(t)\left(r + B^\top p \right) 
- \rho \beta \check{\alpha}^{-1} p 
- \left(1 - \hat{\alpha} \right) \left\lVert A + \left\lvert B \right\rvert \bar{E} \right\rVert_1 \beta p.
\]
Next, using $p(t) \leq \check{\alpha}^{-1} p$ on the right-hand side of the previous inequality, and noting that $1 - \hat{\alpha} = \rho \beta \check{\alpha}^{-1}$, we obtain the desired result:
\[
\check{\alpha}^{-1} \left(1 + \rho \beta \left(1 + \beta \left\lVert A + \left\lvert B \right\rvert \bar{E} \right\rVert_1 \right)\right) p - s 
\ge A^\top p + K^\top(t) \left( r + B^\top p \right).
\]
\jmlrQED
\subsection{Proof of Corollary~\ref{coro1}}\label{Coro1}
\begin{align*}
x^\top(t+1) p &= \left( (A + B K(t)) x(t) + B \epsilon(t) + w(t) \right)^\top p \\
&= x^\top(t) (A + B K(t))^\top p + \left\lvert B \epsilon(t) + w(t) \right\rvert^\top p.
\end{align*}
Define 
\[
\theta \triangleq \check{\alpha}^{-1} \left(1 + \rho \beta \left(1 + \beta \left\lVert A + \left\lvert B \right\rvert \bar{E} \right\rVert_1 \right)\right).
\]
Then, by Theorem~\ref{thm2}, we have
\begin{align*}
(A + B K(t))^\top p 
&\leq \theta p - K^\top(t) r - s \\
&= p - \left(s + K^\top(t) r - (\theta - 1)p \right).
\end{align*}
Let \( \gamma \in (0,1) \) satisfy
\[
s + K^\top(t) r - (\theta - 1)p \ge \gamma \left(s + K^\top(t) r\right),
\]
which holds if
\[
(1 - \gamma)\left(s + K^\top(t) r\right) \ge (\theta - 1)p.
\]
Since \( s + K^\top(t) r \ge s - \bar{E}^\top r \) and \( p \leq \beta s \), it suffices that
\[
(1 - \gamma)\left(s - \bar{E}^\top r\right) \ge (\theta - 1)\beta s,
\]
or equivalently,
\[
\gamma \left(s - \bar{E}^\top r\right) \le s - \bar{E}^\top r - (\theta - 1)\beta s.
\]
Such a \( \gamma \) exists by Assumption~\ref{as:sr}. Therefore,
\[
(A + B K(t))^\top p \leq p - \gamma \left(s + K^\top(t) r\right),
\]
and so
\begin{multline*}
x^\top(t+1) p \leq x^\top(t) p - \gamma x^\top(t)\left(s + K^\top(t) r\right) 
+ \left\lvert B \epsilon(t) + w(t) \right\rvert^\top p \\
\leq x^\top(t) p - \gamma x^\top(t)\left(s + K^\top(t) r\right) 
+ \beta s^\top \left\lvert B \epsilon(t) + w(t) \right\rvert.
\end{multline*}
Rewriting yields
\[
\gamma \left( s^\top x(t) + r^\top K(t) x(t) \right) 
\leq p^\top x(t) - p^\top x(t+1) + \beta s^\top \left\lvert B \epsilon(t) + w(t) \right\rvert.
\]
Summing over \( t_0 \leq t \leq T - 1 \) gives
\[
\sum_{t = t_0}^{T - 1} \left(s^\top x(t) + r^\top K(t) x(t)\right) 
\leq \gamma^{-1} \left( p^\top x_{t_0} 
+ \sum_{t = t_0}^{T - 1} \beta s^\top \left\lvert B \epsilon(t) + w(t) \right\rvert \right).
\]
\jmlrQED

\end{document}